\documentclass[11pt]{amsart}
\usepackage[latin1]{inputenc}
\usepackage{a4}
\usepackage{amssymb}
\usepackage{amsmath}
\begin{document}
\setlength{\parindent}{1.2em}
\def\COMMENT#1{}
\def\TASK#1{}
\def\noproof{{\unskip\nobreak\hfill\penalty50\hskip2em\hbox{}\nobreak\hfill%
       $\square$\parfillskip=0pt\finalhyphendemerits=0\par}\goodbreak}
\def\endproof{\noproof\bigskip}
\newdimen\margin   
\def\textno#1&#2\par{%
   \margin=\hsize
   \advance\margin by -4\parindent
          \setbox1=\hbox{\sl#1}%
   \ifdim\wd1 < \margin
      $$\box1\eqno#2$$%
   \else
      \bigbreak
      \hbox to \hsize{\indent$\vcenter{\advance\hsize by -3\parindent
      \sl\noindent#1}\hfil#2$}%
      \bigbreak
   \fi}
\def\proof{\removelastskip\penalty55\medskip\noindent{\bf Proof. }}
\def\enddiscard{}
\long\def\discard#1\enddiscard{}
\newcommand{\graph}{$K_r^-$}
\newcommand{\fraction}{\frac{r-1}{r(r-2)}}
\newcommand{\comment}{\footnote}
\newcommand{\nat}{\mathbf{N}}
\newtheorem{firstthm}{Proposition}
\newtheorem{thm}[firstthm]{Theorem}
\newtheorem{prop}[firstthm]{Proposition}
\newtheorem{lemma}[firstthm]{Lemma}
\newtheorem{cor}[firstthm]{Corollary}
\newtheorem{problem}[firstthm]{Problem}
\newtheorem{defin}[firstthm]{Definition}
\newtheorem{conj}[firstthm]{Conjecture}
\newtheorem{theorem}[firstthm]{Theorem}
\newtheorem{claim}[firstthm]{Claim}
\def\eps{{\varepsilon}}
\def\N{\mathbb{N}}
\def\R{\mathbb{R}}
\def\K{\mathcal{K}}
\newcommand{\ex}{\mathbb{E}}
\newcommand{\pr}{\mathbb{P}}

\title{Perfect packings with complete graphs minus an edge}
\author{Oliver Cooley \and Daniela K\"uhn \and Deryk Osthus}
\date{}
\maketitle \vspace{-.8cm}
\begin{abstract} \noindent
Let $K_r^-$ denote the graph obtained from $K_r$ by deleting one
edge. We show that for every integer $r\ge 4$ there exists an
integer $n_0=n_0(r)$ such that every graph $G$ whose order $n\ge
n_0$ is divisible by $r$ and whose minimum degree is at least
$(1-1/\chi_{cr}(K_r^-))n$ contains a perfect $K_r^-$-packing, i.e.~a
collection of disjoint copies of $K_r^-$ which covers all vertices
of~$G$. Here $\chi_{cr}(K_r^-)=\frac{r(r-2)}{r-1}$ is the critical
chromatic number of~$K_r^-$. The bound on the minimum degree is best
possible and confirms a conjecture of Kawarabayashi for large~$n$.
\end{abstract}

\section{Introduction}\label{intro}
Given two graphs $H$ and $G$, an \emph{$H$-packing in $G$} is a
collection of vertex-disjoint copies of $H$ in $G$. An $H$-packing
in $G$ is called \emph{perfect} if it covers all vertices of $G$. In
this case, we also say that $G$ contains an \emph{$H$-factor}. The
aim now is to find natural conditions on $G$ which guarantee the
existence of a perfect $H$-packing in $G$. For example, a famous
theorem of Hajnal and Szemer\'edi~\cite{HSz70} gives a best possible
condition on the minimum degree of $G$ which ensures that $G$ has a
perfect $K_r$-packing. More precisely, it states that every graph
$G$ whose order $n$ is divisible by $r$ and whose minimum degree
is at least $(1-1/r)n$ contains a perfect $K_r$-packing. (The
case $r=3$ was proved earlier by Corr\'adi and Hajnal~\cite{CH}
and the case $r=2$ follows immediately from Dirac's theorem on
Hamilton cycles.)

Alon and Yuster~\cite{AY96} proved an extension of this result to
perfect packings of arbitrary graphs~$H$. They showed that for every
$\gamma>0$ and each graph $H$ there exists an integer
$n_0=n_0(\gamma,H)$ such that every graph $G$ whose order $n\ge n_0$
is divisible by $|H|$ and whose minimum degree is at least
$(1-1/\chi(H)+\gamma)n$ contains a perfect $H$-packing. They
observed that there are graphs $H$ for which the error term $\gamma
n$ cannot be omitted completely, but conjectured that it could be
replaced by a constant which depends only on~$H$. This conjecture
was proved by Koml\'os, S\'ark\"ozy and Szemer\'edi~\cite{KSSz01}.

Thus one might think that just as in Tur\'an theory -- where instead
of an $H$-packing one only asks for a single copy of~$H$ -- the
chromatic number of $H$ is the crucial parameter when one considers
$H$-packings. However, one indication that this is not the case is
provided by the result of Koml\'os~\cite{JKtiling}, which states
that if one only requires an \emph{almost} perfect $H$-packing
(i.e.~one which covers almost all of the vertices of $G$), then the
relevant parameter is the criticial chromatic number of~$H$. Here
the \emph{critical chromatic number} $\chi_{cr}(H)$ of a graph $H$
is defined as $(\chi(H)-1)|H|/(|H|-\sigma(H))$, where $\sigma(H)$
denotes the minimum size of the smallest colour class in a colouring
of $H$ with $\chi(H)$ colours and where $|H|$ denotes the order of
$H$. Note that $\chi_{cr}(H)$ always satisfies  $\chi(H)-1 <
\chi_{cr}(H) \le \chi(H)$ and is closer to $\chi(H)-1$ if
$\sigma(H)$ is comparatively small. Building on this, in~\cite{HPack}
it was shown that for some graphs~$H$ the critical chromatic
number is even the relevant parameter for perfect packings, while
for all other graphs the relevant parameter is the chromatic number.
In order to state the precise result (Theorem~\ref{the:hpack}) we
need to introduce some notation. A colouring of a graph $H$ is
called \textit{optimal} if it uses exactly $\chi(H)$ colours. Let
$\ell := \chi(H)$. Given an optimal colouring~$c$ of $H$, let $x_1
\leq x_2 \leq \ldots \leq x_\ell$ be the sizes of the colour
classes. Define $\mathcal{D}(c)=\lbrace x_{i+1}-x_i \mid
i=1,\ldots,\ell-1 \rbrace$. Let $\mathcal{D}(H)$ be the union of all
the sets $\mathcal{D}(c)$ over all optimal colourings~$c$ of~$H$. We
define $hcf_\chi(H)$ to be the highest common factor of the elements
of~$\mathcal{D}(H)$ (or $hcf_\chi(H):=\infty$ if $\mathcal{D}(H) =
\lbrace 0 \rbrace$). Define $hcf_c(H)$ to be the highest common
factor of the orders of all the components of $H$. For any graph
$H$, if $\chi(H) \neq 2$, we say $hcf(H)=1$ if $hcf_\chi(H)=1$. If
$\chi(H)=2$, we say $hcf(H)=1$ if both $hcf_c(H)=1$ and
$hcf_\chi(H)\leq 2$.

\begin{thm} \label{the:hpack}{\bf \cite{HPack}}
Given a graph $H$, let $\delta(H,n)$ denote the smallest integer $k$
such that every graph $G$ whose order $n$ is divisible by~$|H|$ and
with $\delta(G)\ge k$ contains a perfect $H$-packing. Then
$$\delta(H,n) =
\begin{cases}
 \left(1-\frac{1}{\chi_{cr}(H)} \right)n+O(1) &\text{if hcf$(H)=1$}, \\
\left(1-\frac{1}{\chi(H)} \right)n+O(1) &\text{if hcf$(H)\neq 1$}. \\
\end{cases} $$
\end{thm}
Here the $O(1)$ error term depends only on $H$ and there are graphs
$H$ for which it cannot be omitted completely (see
Proposition~\ref{properrorterm}). Also, note that the upper bound on
$\delta(H,n)$ in the case when $hcf(H) \neq 1$ is the result
in~\cite{KSSz01} mentioned earlier. The proof in~\cite{HPack}
for the case when $hcf(H) = 1$ gave a
constant which was dependent on the constant in Szemer\'{e}di's
regularity lemma, and is therefore huge.

Our main result shows that in the case when $H=K_r^-$, where $r \geq
4$, the error term in Theorem~\ref{the:hpack} can be omitted
completely. (Recall that $K_r^-$ denotes the graph obtained from
$K_r$ by deleting one edge.) Note that $hcf(K_r^-)=1$ for $r \geq
4$.

\begin{thm}\label{Klminusthm}\label{the:main}
For every integer $r\ge 4$ there exists an integer $n_0=n_0(r)$ such
that every graph $G$ whose order $n\ge n_0$ is divisible by $r$ and
whose minimum degree is at least
$$\left( 1-\frac{1}{\chi_{cr}(K_r^-)}\right)n$$ contains a perfect
$K_r^-$-packing.
\end{thm}
This theorem confirms a conjecture of Kawarabayashi~\cite{KK} for large~$n$.
The case $r=4$ of
the conjecture (and thus of Theorem~\ref{Klminusthm}) was proved by
Kawarabayashi~\cite{KK}. By a result of Enomoto, Kaneko and
Tuza~\cite{EKT87}, the conjecture also holds for the case $r=3$
under the additional assumption that $G$ is connected. (Note that
$K_3^-$ is just a path on $3$ vertices and that in this case the
required minimum degree equals $n/3$.)
For completeness, in Proposition~\ref{exexample} we will
give an explicit construction
showing that the bound on the minimum degree in Theorem~\ref{Klminusthm}
is best possible.

Clearly, it would be desirable to characterize all those graphs for
which the $O(1)$-error term in Theorem~\ref{the:hpack} can be
omitted.
However, we do not know what such a characterization might look like.
By the Hajnal-Szemer\'edi theorem~\cite{HSz70} the error term can be
omitted for complete graphs. A result of Abbasi~\cite{Abbasi}
implies that, for large~$n$, it can be omitted for cycles. In~\cite{MPhil} the
first author describes a further class of graphs for which
the ideas in this paper can be adapted to remove the error term
completely for large~$n$. On the other hand,
Proposition~\ref{properrorterm} shows that the error term cannot be
omitted if $H$ is a complete $\ell$-partite graph with $\ell \geq 3$
and at least $\ell-1$ vertex classes of size at least~3. A larger class
of graphs~$H$ for which this is the case is given in~\cite{MPhil}.

Algorithmic issues related to Theorem~\ref{the:hpack} are discussed
in~\cite{KOSODA}. It was shown there that for any $\varepsilon
>0$ the perfect $H$-packing guaranteed by Theorem~\ref{the:hpack}
can be found in polynomial time if the $O(1)$-error term is replaced by
$\varepsilon n$. Moreover,
if the minimum degree condition on $G$ is reduced
a little below the threshold, then there are many graphs $H$ for
which the decision problem of whether $G$ has a perfect $H$-packing becomes
NP-complete.

\section{Notation and preliminaries}\label{sec:prelim}
Throughout this paper we omit floors and ceilings whenever this does
not affect the argument. We write $e(G)$ for the number of edges of
a graph $G$, $|G|$ for its order, $\delta(G)$ for its minimum
degree, $\Delta(G)$ for its maximum degree, $\chi(G)$ for its
chromatic number and $\chi_{cr}(G)$ for its critical chromatic
number as defined in Section~\ref{intro}. We denote the degree of a
vertex $x\in G$ by $d_G(x)$ and its neighbourhood by~$N_G(x)$. Given
a vertex set $A\subseteq V(G)$, we also write $N_A(x)$ for the set
of all neighbours of $x$ in~$A$. We denote by $G[A]$ the subgraph of
$G$ induced by the vertex set $A$. Given disjoint sets $A,B\subseteq
V(G)$, we denote by $e(A,B)$ the number of all edges between $A$ and
$B$ and write $d(A,B):=e(A,B)/|A||B|$ for the density of the
bipartite subgraph of~$G$ between~$A$ and~$B$. We denote by
$d(A):=e(A)/\binom{|A|}{2}$ the density of~$A$.

%
For a graph $H$ of chromatic number~$\ell$, define
the \textit{bottle graph $B^*(H)$ of $H$}, to be the complete
$\ell$-partite graph which has $\ell-1$ classes of size
$|H|-\sigma(H)$ and one class of size $(\ell-1)\sigma(H)$. (Recall
that $\sigma(H)$ is the smallest possible size of a colour class in
an $\ell$-colouring of $H$.) Thus $B^*(H)$ contains a perfect $H$-packing
consisting of $\ell-1$ copies of~$H$.
We will use $B^*$ to denote $B^*(K_r^-)$ whenever this is unambiguous.

For completeness, we include the construction which shows that the
bound on the minimum degree in Theorem~\ref{Klminusthm} is best
possible.

\begin{prop}\label{exexample}
Let $r \ge 4$. Then for all $k \in \mathbb{N}$ there is a graph $G$
on $n=kr$ vertices whose minimum degree is
$\left\lceil(1-1/\chi_{cr}(K^-_r))n\right\rceil-1$ but which does
not contain a perfect $K_r^-$-packing.
\end{prop}
\proof We construct $G$ as follows. $G$ is a complete
$(r-1)$-partite graph with vertex classes $U_0,\dots,U_{r-2}$, where
$|U_0|=k-1$ and the sizes of all other classes are as equal as
possible. It is easy to check that $G$ has the required
minimum degree.%
     \COMMENT{Indeed,
$$\delta(G) =n-\left\lceil \frac{n-|U_0|}{r-2} \right\rceil
=k(r-1)-\left\lceil\frac{k+1}{r-2}\right\rceil
\ge k(r-1)-\left(\frac{k}{r-2}+1\right)
= \frac{r^2-3r+1}{r(r-2)}n-1.$$
}
Moreover, every copy of $K_r^-$ in $G$ contains at least one vertex
in $U_0$. Thus we can find at most $|U_0|$ pairwise disjoint copies
of $K_r^-$ which therefore cover at most $(k-1)(r-1)<n-|U_0|$
vertices of $G-U_0$. Thus $G$ does not contain a perfect
$K_r^-$-packing.
\endproof

Note that Proposition~\ref{exexample} extends to every graph~$H$
which is obtained from a~$K_{r-1}$ by adding a new vertex and
joining it to at most $r-2$ vertices of the~$K_{r-1}$. Since each
such~$H$ is a subgraph of $K_r^-$ and since
$\chi_{cr}(H)=\chi_{cr}(K_r^-)$, it follows from this observation
and from Theorem~\ref{the:main} that
$\delta(H,n)=\lceil(1-1/\chi_{cr}(H))n\rceil$ if~$n$ is sufficiently
large (where $\delta(H,n)$ is as defined in
Theorem~\ref{the:hpack}).

The following example shows that for a large class of graphs, the
$O(1)$-error term in Theorem~\ref{the:hpack} cannot be omitted
completely. The example is an extension of a similar construction
in~\cite{KSSz01}.

\begin{prop}\label{properrorterm}
Suppose that $H$ is a complete $\ell$-partite graph with $\ell \geq 3$
such that every vertex class of~$H$, except possibly its smallest class, has
at least $3$ vertices. Then there are infinitely many graphs
$G$ whose order~$n$ is divisible by~$|H|$, whose minimum degree
satisfies $\delta(G)=(1-\frac{1}{\chi_{cr}(H)})n$ but which do not
contain a perfect $H$-packing.
\end{prop}
\proof
Let $\sigma$ denote the size of the smallest vertex class of~$H$.
Given $k\in\N$, consider the
complete $\ell$-partite graph on $n:= k(\ell-1)|H|$ vertices
whose vertex classes $A_1,\dots,A_\ell$ satisfy $|A_1|:=(|H|-\sigma)k+1$,
$|A_\ell|:= k(\ell-1)\sigma-1$ and $|A_i|:=(|H|-\sigma)k$
for all $1<i<\ell$. Let $G$ be the graph obtained by adding
a perfect matching into~$A_1$ or, if $|A_1|$ is odd, a matching covering
all but 3 vertices and a path of length~2 on these remaining vertices.
Observe that the minimum degree of $G$ is $(1-\frac{1}{\chi_{cr}(H)})n$.%
      \COMMENT{$\delta(G)=(\ell-2)(|H|-\sigma)k + \sigma k (\ell-1)$\\\\
$\frac{1}{\chi_{cr}(H)}= \frac{|H|-\sigma}{(\ell-1)|H|}$\\\\
$n=|H|k(\ell - 1)$\\\\
$(1-\frac{1}{\chi_{cr}(H)})n = (\ell -1)|H|k - (|H|-\sigma)k =
(\ell -2)|H|k + \sigma k = \delta(G)$.}

Consider any copy~$H'$ of~$H$ in~$G$. Suppose that~$H'$ meets~$A_\ell$
in at most $\sigma - 1$ vertices. Then there is a colour class~$X$ of~$H'$
which meets $A_\ell$ but does not lie entirely in~$A_\ell$.
So some vertex class of~$G$ must meet at least two colour classes
of~$H'$. Since $H'$ is complete $\ell$-partite, this vertex class
must have some edges in it, and so must be~$A_1$. However, $A_1$
cannot meet three colour classes of~$H'$, since it is triangle free.
Thus every colour class of~$H'$ except $X$ lies completely within
one $A_i$. Furthermore, $A_1$ cannot contain two complete colour
classes of~$H'$, since then $G[A_1]$ would have a vertex of degree
$3$, a contradiction. So $A_1$ meets $X$ as well as another colour
class~$Y$ of~$H'$. Furthermore $X\setminus A_\ell\subseteq A_1$
and $Y\subseteq A_1$.%
     \COMMENT{If $X\setminus A_\ell\not \subseteq A_1$ then the 3rd vertex class
of $G$ which meets $X$ avoids all the other colour classes of $H'$
and thus $A_1$ would meet $X$, $Y$ and another colour class of $H'$,
which is impossible since $A_1$ doesn't contain a triangle. A
similar argument works for~$Y$.} Let $x\in X\cap A_1$. Then
$Y\subseteq N_G(x)$ since $Y \subseteq N_{H'}(x)$. This implies that
$|Y|\le 2$ and so $\sigma=|Y|\le 2$. Thus $|X|\ge 3$. Since at most
$\sigma-1\le 1$ vertices of~$X$ lie in~$A_\ell$ this in turn implies
that $|X\cap A_1|\ge 2$. As $X\cap A_1$ lies in the neighbourhood of
any vertex from~$Y$, we must have that $|X\cap A_1|=2$. Thus $X\cap
A_1$ can only lie in the neighbourhood of one vertex from~$Y$. Hence
$\sigma = |Y| = 1$. But then $X$ avoids $A_\ell$, a contradiction.

So any copy of~$H$ in~$G$ has at least $\sigma$ vertices in~$A_\ell$.
Thus any $H$-packing in~$G$ consists of less than $k(\ell-1)$
copies of~$H$ and therefore covers less than $k(\ell-1)(|H|-\sigma)<|G|-|A_\ell|$
vertices of $G-A_\ell$.
So $G$ does not contain a perfect $H$-packing.
\endproof

Note that the proof of Proposition~\ref{properrorterm}
shows that if $|H|-\sigma$ is odd then we only need that
every vertex class of~$H$ (except possibly its smallest class)
has at least two vertices.%
     \COMMENT{We then have to take an odd integer $k$.}
Moreover, it is not hard to see that the conclusion of Proposition~\ref{properrorterm} holds
for all graphs~$H$ which do not have an optimal colouring with a vertex class
of size~$\sigma+1$ (see~\cite{MPhil} for details).

In the proof of Theorem~\ref{Klminusthm} we will use the following
observation about packings in almost complete $(q+1)$-partite
graphs. It follows easily from the Blow-up lemma (see e.g.~\cite{JKblowup}),
but we also sketch how it can be deduced directly from Hall's
theorem.

\begin{prop}\label{almostcomplete}
For all $q,r\in \N$ there exists a positive constant
$\tau_0=\tau_0(q,r)$ such that the following holds for every
$\tau\le \tau_0$ and all $k\in\N$. Let $H_{q,r}$ be the complete
$(q+1)$-partite graph with $q$ vertex classes of size $r$ and one
vertex class of size~$1$. Let $G^*$ be a $(q+1)$-partite graph with
vertex classes $V_1,\dots,V_{q+1}$ such that $|V_i|=kr$ for all
$i\le q$ and such that $|V_{q+1}|=k$. Suppose that for all distinct
$i,j\le q+1$ every vertex $x\in V_i$ of $G^*$ is adjacent to all but
at most $\tau |V_j|$ vertices in $V_j$. Then $G^*$ has a perfect
$H_{q,r}$-packing.
\end{prop}
\proof We proceed by induction on~$q$. If $q=1$ then we are looking
for a perfect $K_{1,r}$-packing. So the result can easily be deduced
from Hall's theorem with $\tau_0=1/2$.%
     \COMMENT{$|V_1|=kr$,
$|V_2|=k$. Blow up $V_2$ by a factor of $r$ and look for a perfect
matching, so check Hall's condition.
Suppose Hall's condition fails for $A \subseteq V_1$.
Then since $\tau_0 = \frac{1}{2}$, $|N(A)|\geq \frac{1}{2}kr$, so $|A|>\frac{1}{2}kr$.
Now $\emptyset \neq N(V_2 \backslash N(A)) \subseteq V_1 \backslash A$
So $|V_1 \backslash A| \geq |N(V_2 \backslash N(A))| \geq
\frac{1}{2}kr$
So $|A| \le \frac{1}{2}kr$ - contradiction}
Now suppose that $q>1$ and let
$\tau_0(q,r)\ll\tau_0(q-1,r)$. As before, we can find a perfect
$K_{1,r}$-packing in $G^*[V_q\cup V_{q+1}]$. Let $G'$ be the graph
obtained from $G^*$ by replacing each copy $K$ of such a $K_{1,r}$
with one vertex $x_K$ and joining $x_K$ to $y\in V_1\cup\dots\cup V_{q-1}$
whenever $y$ is adjacent to every vertex of~$K$. Then $G'$ contains
a perfect $H_{q-1,r}$-packing by induction. Clearly, this
corresponds to a perfect $H_{q,r}$-packing in~$G^*$.
\endproof
\section{Overview of the proof}\label{sec:overview}
Our main tool is the following result from~\cite{HPack}. It states that
in the ``non-extremal case'', where the graph~$G$ given in Theorem~\ref{the:hpack}
satisfies certain conditions, we can find a perfect packing even if the
minimum degree is slightly smaller than required in Theorem~\ref{the:hpack}.
The conditions ensure that the graph $G$ does not
look too much like one of the extremal examples of graphs whose
minimum degree is just a little smaller than required in Theorem~\ref{the:hpack}
but which do not contain a perfect $H$-packing.

\begin{theorem}
\label{the:nonex}
Let $H$ be a graph of chromatic number $\ell \geq 2$ with $hcf(H)=1$.
Let $z_1$ denote the size of the small class of the bottle graph
$B^*(H)$, let $z$ denote the size of one of the large classes, and let
$\xi = z_1/z$. Let $\theta \ll \tau_0 \ll \xi,1-\xi,1/|B^*(H)|$
be positive constants. There exists an integer $n_0$ such that the
following holds. Suppose $G$ is a graph whose order $n \geq n_0$ is
divisible by $|B^*(H)|$ and whose minimum degree satisfies
$\delta(G)\geq (1-\frac{1}{\chi_{cr}(H)}-\theta)n$. Suppose that $G$
also satisfies the following conditions:
\begin{itemize}
\item[{\rm (i)}] $G$ does not contain a vertex set $A$ of size $zn/|B^*(H)|$
such that $d(A)\leq \tau_0$.
\item[{\rm (ii)}]  If $\ell = 2$, then $G$ does not contain a vertex set $A$ with
$d(A,V(G)\setminus A) \leq \tau_0$.
\end{itemize}
Then $G$ has a perfect $H$-packing.
\end{theorem}
By applying this theorem with $H:=K_r^-$ (where $r\ge 4$), we only
need to consider the extremal case, when there are large almost
independent sets.
(Note that if the order of the graph $G$ given by Theorem~\ref{the:main} is not
divisible by $|B^*(K_r^-)|$, we must first greedily
remove some copies of $K_r^-$ before applying
Theorem~\ref{the:nonex}. The existence of these copies follows from
the Erd\H{o}s-Stone theorem, and since we only need to remove a
bounded number of copies, this will not affect any of the properties
required in Theorem~\ref{the:nonex} significantly.)

Suppose that we have $q$ such large almost independent sets. Then we
will think of the remainder of the vertices of $G$ as the $(q+1)$th
set. We will show in Section \ref{sec:tidy} that by taking out a few
copies of $K_r^-$ and rearranging these $q+1$ sets slightly, we can
achieve that these sets will induce an almost complete
$(q+1)$-partite graph. Furthermore, the proportion of the size of
each of the first $q$ of these modified sets to the size of the
entire graph will be the same as for the large classes of the bottle
graph $B^*(K_r^-)$ defined in Section~\ref{sec:prelim}.

Let $B_1^*$ be the subgraph of $B^*(K_r^-)$ obtained by deleting~$q$ of the
large vertex classes. Ideally, we would like to apply Theorem~\ref{the:nonex} to
find a $B_1^*$-packing in the (remaining) subgraph of~$G$ induced by the
$(q+1)$th vertex set. In a second step we would then like to extend this $B_1^*$-packing
to a $B^*(K_r^-)$-packing in~$G$, using the fact that the $(q+1)$-partite
subgraph of $G$ between the classes defined above is almost complete. This would
clearly yield a \graph-packing of~$G$.

However, there are some difficulties. For example,
Theorem~\ref{the:nonex} only applies to graphs~$H$ with
$hcf(H)=1$, and this may not be the case for~$B_1^*$ if it is
bipartite. So instead of working with $B_1^*$, we consider a
suitable subgraph~$B_1$ of~$B^*_1$ which does satisfy $hcf(B_1)=1$.
Moreover, if $B_1$ is bipartite we may have to take out a few
further carefully chosen copies of $K_r^-$ from $G$ to ensure that
condition~(ii) is also satisfied before we can apply
Theorem~\ref{the:nonex} to the subgraph induced by the
$(q+1)$th vertex set.

\section{Tidying up the classes}
\label{sec:tidy} Let $n$ and $q$ be integers such that $n$ is
divisible by $r(r-2)=|B^*(K_r^-)|$ and such that $1\le q \le r-2$.
Note that in the case when $H:=K_r^-$ the set~$A$ in condition~(i)
of Theorem~\ref{the:nonex} has size $\frac{r-1}{r(r-2)}n$.
We say that disjoint vertex sets $A_1,\dots,A_{q+1}$ are
\emph{$(q,n)$-canonical} if $|A_i|=\frac{r-1}{r(r-2)}n$ for all
$i\le q$ and $|A_{q+1}|=\frac{n}{r} + (r-q-2)\frac{r-1}{r(r-2)}n
=n-\sum_{i=1}^q|A_i|$. Note that in this case the graph $K(q,n)$
obtained from the complete graph on $\bigcup_{i=1}^{q+1} A_i$ by
making each $A_i$ with $i\le q$ into an independent set has a
perfect $B^*(K_r^-)$-packing and thus also a perfect
$K_r^-$-packing.

Our aim in the
following lemma is to remove a few disjoint copies of $K_r^-$ from our given
graph $G$ in order to obtain a graph on $n^*$ vertices which looks almost like
$K(q,n^*)$. In the next section we will then use this property to show that
this subgraph of~$G$ has a perfect $K_r^-$-packing.

\begin{lemma}
\label{lem:tidy}
Let $r \geq 4$ and $0<\tau \ll 1/r$. Then there exists an integer
$n_0=n_0(r,\tau)$ such that the following is true.
Let $G$ be a graph whose order $n\ge n_0$ is divisible by~$r$ and whose
minimum degree satisfies $\delta(G) \geq (1-\frac{1}{\chi_{cr}(K_r^-)})n$.
Suppose that for some $1 \leq q \leq r-2$ there are $q$ disjoint vertex sets
$A_1,\ldots,A_q$ in $G$ such that $|A_i|=\lceil \fraction n \rceil$ and
$d(A_i)\leq \tau$ for $1 \leq i \leq q$. Set
$A_{q+1}:=V(G)\backslash (A_1 \cup \ldots \cup A_q)$. Then there exist
disjoint vertex sets $A_1^*,\ldots,A_{q+1}^*$ such that the following hold:
\begin{itemize}
\item[{\rm (i)}] If $G^*:= G[\bigcup_{i=1}^{q+1}A_i^*]$ and $n^*:=|G^*|$ then $r(r-2)$
divides $n^*$, and $G-G^*$ contains a perfect \graph-packing. Furthermore,
$n-n^*\leq \tau^{1/3}n$.
\item[{\rm (ii)}] $|A_1^*|=|A_2^*|=\ldots=|A_q^*|=\fraction n^*$.
\item[{\rm (iii)}] For all $i,j \leq q+1$ with $i \neq j$, each vertex in $A_i^*$
has at least $(1-\tau^{1/5})|A_j^*|$ neighbours in $A_j^*$.
\end{itemize}
\end{lemma}
\proof Note that if $n$ is divisible by $r(r-2)$ then the sets
$A_1,\dots,A_{q+1}$ are $(q,n)$-canonical. If $n$ is not divisible
by $r(r-2)$ then we will change the sizes of the $A_i$ slightly as
follows. Write $n=n'+kr$ where $n'$ is divisible by $r(r-2)$ and
$0<k < r-2$. If $k\ge q$ then we do not change the sizes of
the~$A_i$. If $k<q$ then for each $i$ with $k<i\le q$ we move one
vertex from $A_i$ to $A_{q+1}$. We still denote the sets thus
obtained by $A_1,\dots,A_{q+1}$. We may choose the vertices we move
in such a way that the density of each $A_i$ with $i\le q$ is still
at most~$\tau$. Note that $\lceil \fraction n \rceil = \fraction
n'+k+1$. Thus both in the case when $k\ge q$ and in the case when
$k<q$ the sets $A_1,\dots,A_{q+1}$ can be obtained from
$(q,n')$-canonical sets by adding $kr$ new vertices as follows. For
each $i\le \min\{k,q\}$ we add $k+1$ of the new vertices to the
$i$th vertex set, for each $i$ with $\min\{k,q\}<i\le q$ we add $k$
new vertices to the $i$th vertex set and all the remaining new
vertices are added to~$A_{q+1}$. Let $K$ be the graph obtained from
the complete graph on $\bigcup_{i=1}^{q+1} A_i$ by making each $A_i$
with $i\le q$ into an independent set. It is easy to see that
$K(q,n')$ can be obtained from $K$ by removing $k$ disjoint copies
of~$K_r^-$. In particular, $K$ has a perfect $K_r^-$-packing. Note
that if $k<q$ then this would not hold if we had not changed the
sizes of the~$A_i$. Later on we will use that in all cases we have
\begin{equation}\label{eqsizeAi}
|A_i| \ge \fraction n'+k=\fraction (n-kr)+k
\end{equation}
for all $i\le q$, where we set $n':=n$ and $k:=0$ if $n$ is divisible by $r(r-2)$.
Observe that $\chi_{cr}(K_r^-)=\frac{r(r-2)}{r-1}$ and so
$\delta(G)\ge (1-\frac{r-1}{r(r-2)})n$.
Thus the minimum degree condition on $G$ implies that the neighbours of
any vertex might essentially avoid one of the $A_i$, for $i\leq q$, but no more.

Now for each index $i$, call a vertex $x \in A_i$ \textit{i-bad} if
$x$ has at least $\tau^{1/3}|A_i|$ neighbours in $A_i$. Note that,
for $i \leq q$, the number of $i$-bad vertices is at most
$\tau^{2/3}|A_i|$ since $d(A_i)\le \tau$ for such~$i$. Call a vertex
$x \in A_i$ \textit{i-useless} if, for some $j \neq i$, $x$ has at
most $(1-\tau^{1/4})|A_j|$ neighbours in $A_j$. In this case the
minimum degree condition shows that, provided $i \neq r-1$, $x$ must
have at least a $\tau^{1/3}$-fraction of the vertices in its own
class as neighbours, i.e. $x$ is $i$-bad. Thus every vertex that is
$i$-useless is also $i$-bad for $i \neq r-1$. In particular, for
each $i \leq q$, there are at most $\tau^{2/3}|A_i|$ $i$-useless
vertices.

For $i=q+1$ we estimate the number $u_{q+1}$ of $(q+1)$-useless vertices by
looking at the edges between $A_{q+1}$ and $V(G)\backslash A_{q+1}$.
We have
\begin{align*}
e(A_{q+1},V(G)\backslash A_{q+1}) & \geq \sum_{i=1}^{q}\lbrace
|A_i| \delta(G) - 2e(A_i) - \sum_{j \neq i, j\le q}|A_i||A_j|  \rbrace \\
& \geq q(|A_1|-1)\delta(G) -q\tau|A_1|^2 - q(q-1)|A_1|^2.
\end{align*}
On the other hand,
\begin{align*}
e(A_{q+1},V(G)\backslash A_{q+1}) & \leq u_{q+1}\lbrace(q-1)|A_1|+
(1-\tau^{1/4})|A_1|\rbrace + (|A_{q+1}|-u_{q+1})q|A_1|\\
& = q|A_1||A_{q+1}| - u_{q+1}\tau^{1/4}|A_1|.
\end{align*}
Combining these inequalities gives, after some calculations,%
     \COMMENT{$\tau^{\frac{1}{4}}u_{q+1} \leq q|A_{q+1}| + q\tau |A_1|
     - \frac{1}{|A_1|} [ q(|A_1|-1)\delta(G)-q(q-1)|A_1|^2]\\
     \\
     \delta(G)-(q-1)|A_1| \geq |A_{q+1}| - q\\$
     \begin{align*}
     \tau^{\frac{1}{4}}u_{q+1} & \leq q|A_{q+1}| + q\tau|A_1| -
     q(|A_{q+1}|-q) + q\frac{\delta(G)}{|A_1|}\\
     & = q^2 + q\tau |A_1| + q\frac{\delta(G)}{|A_1|}\\
     & \leq q^2 + q\tau |A_1| + rq\\
     & \leq 3q \tau |A_1|\\
     u_{q+1} & \leq \tau^\frac{2}{3} |A_{q+1}|
     \end{align*}
    }
that $u_{q+1}\leq \tau^{2/3}|A_{q+1}|$. So in total the number of
vertices which are $i$-useless for some $i$ is at most
$\tau^{2/3}n$.

Given $j \neq i$, call a vertex $x \in A_i$ \textit{j-exceptional}
if $x$ has at most $\tau^{1/3}|A_j|$ neighbours in $A_j$. Thus every
such vertex is also $i$-useless, and therefore $i$-bad if $i < r-1$.
Furthermore, if $i=r-1$, then an exceptional vertex in $A_i$ is also
$i$-bad. So all exceptional vertices are bad.

Now if for some $i\neq j$ there exists an $i$-bad vertex $x \in A_i$
and an $i$-exceptional vertex $y \in A_j$, then let us swap $x$ and
$y$. (Note that a vertex is not $i$-exceptional for more than
one~$i$.) Having done this, since there are not too many exceptional
vertices, we will still have that each non-bad vertex in $A_i$ has
at most $2\tau^{1/3}|A_i|$ neighbours in $A_i$, each non-useless
vertex in $A_i$ still has at least $(1-2\tau^{1/4})|A_j|$ neighbours
in each $A_j$ with $j \neq i$ and each non-$i$-exceptional vertex
still has at least $\tau^{1/3}|A_i|/2$ neighbours in $A_i$. We will
also have that for any $i$ for which $i$-exceptional vertices exist,
there are no $i$-bad vertices.

We now wish to remove all the exceptional vertices by taking out a few
disjoint copies of $K_r^-$ which will cover them. For simplicity, we
will split the argument into two cases. In both cases we will repeatedly
remove $r-2$ disjoint copies of $K_r^-$ at a time. We say that such a collection
of $r-2$ copies \emph{respects the proportions of the~$A_i$} if
altogether these copies meet each $A_i$ with $i\le q$ in exactly $r-1$ vertices.

\medskip

\noindent
\textbf{Case 1.} $q \leq r-3$

\smallskip

\noindent In this case the minimum degree condition ensures that no
vertex is $(q+1)$-exceptional. To deal with the $j$-exceptional vertices
for $j \leq q$ we will need the fact that we
can find a reasonably large number of disjoint copies of $K_{r-1-q}$
in $G[A_{q+1}]$. To prove this fact, observe that

\begin{equation}\label{eqmindeg1}
\delta(G[A_{q+1}])  \geq  \delta(G)-\sum_{i=1}^q |A_i|
 \geq  |A_{q+1}| - \fraction n
\end{equation}
and
\begin{equation}\label{eqmindeg2}
\fraction \frac{n}{|A_{q+1}|} \stackrel{(\ref{eqsizeAi})}{\leq}
\fraction \frac{1}{\frac{1}{r}+(r-q-2) \fraction} \leq
\frac{1}{r-q-2}-c(r)
\end{equation}
where $c(r)>0$ is a constant depending only on $r$. Combining these
results gives
\begin{equation}\label{eqmindegAq+1}
\delta(G[A_{q+1}]) \geq
\left(1-\frac{1}{r-q-2}+c(r)\right)|A_{q+1}|.
\end{equation}
Thus we can apply Tur\'an's theorem repeatedly to find
at least $\frac{c(r)}{r-q-1}|A_{q+1}|$ disjoint copies of
$K_{r-q-1}$ in~$G[A_{q+1}]$.

Now for each $i\le q+1$ in turn, consider the exceptional vertices
$x\in A_i$. Suppose that~$x$ is $j$-exceptional. First move~$x$ into~$A_j$.
Note that the minimum degree condition
on $G$ means that $x$ is joined to almost all vertices in~$A_\ell$
for every $\ell\neq j$. We greedily choose a copy of $K_r^-$
covering $x$ and one other vertex in $A_j$, $r-q-1$ vertices in
$A_{q+1}$ and one vertex in all other classes, where all vertices
other than $x$ were chosen to be non-useless. (Indeed, to find such a
copy of $K_r^-$ we first choose a copy of $K_{r-q-1}$ in $A_{q+1}$
which lies in the neighbourhood of~$x$ and which consists of non-useless
vertices. Then we choose all the remaining vertices.)%
     \COMMENT{If $i\le q$ we remove 2 vs from $A_i$. If $i=q+1$ we remove
$r-q$ vertices from~$A_{q+1}$. This is also ok since $q\le r-3$.}
Remove this copy of \graph. Also greedily remove $r-3$
further disjoint copies of $K_r^-$ such that together all these
copies of $K_r^-$ respect the proportions of the~$A_i$.
Proceed similarly for all the exceptional vertices.
For each exceptional vertex we are removing $r-2$ copies of \graph,
so in total we are removing at most $r(r-2)\tau^{2/3}n$ vertices.

\medskip

\noindent
\textbf{Case 2.} $q=r-2$

\smallskip

\noindent
In this case, the exceptional vertices in $A_{r-1}$ need special attention
since we cannot simply move them into another class without making $A_{r-1}$
too small. So we proceed as follows.
For each $i \leq r-2$, let $s_i$ be the number of $i$-exceptional
vertices in $A_{r-1}$. Whenever $s_i>0$ we will find a matching of size~$s_i$
in $G[A_i]$. To see that such a matching exists,
consider a maximal matching in $A_i$ and let~$m$ denote the size of this
matching. Note that
\[e(A_i) \leq  2m \Delta(A_i) \leq 2m 2\tau^{1/3}|A_i|\]
since the presence of $i$-exceptional vertices guarantees that no
vertex in~$A_i$ is $i$-bad. Also
\begin{align*}
e(A_i) & \geq \frac{1}{2}\lbrace \delta(G)|A_i|-(n-|A_i|-s_i)|A_i|
- s_i 2\tau^{1/3}|A_i|\}\\
& \ge \frac{|A_i|}{2}\lbrace |A_i| - \fraction n +
s_i(1-2\tau^{1/3})\}\\
& \stackrel{(\ref{eqsizeAi})}{\geq} \frac{|A_i|}{2} \lbrace s_i(1-2\tau^{1/3})-
\frac{k}{r-2}\}.
\end{align*}
Since $k \leq r-3$ and $\tau \ll 1/r$, comparing these two bounds on
$e(A_i)$ gives $m \gg s_i$ whenever $s_i > 0$. So we may pick a
matching $M_i$ with $s_i$ edges in $A_i$, all of whose vertices are
non-useless (since no vertices in $A_i$ are bad). Now for each $i$
in turn, we will remove the $i$-exceptional vertices in $A_{r-1}$
using this matching. For each such vertex $x\in A_{r-1}$, pick an
edge $yz\in M_i$. Swap $x$ with $y$; we now no longer consider $x$
to be exceptional. Then greedily find a copy of $K_r^-$ which meets
$A_{r-1}$ precisely in~$y$, which meets $A_i$ precisely in~$z$ and
which contains two vertices in some $A_j$ with $j\neq i,r-1$ (such a
$j$ exists since $r \geq 4$), and one vertex in each other $A_j$.
All these vertices will be chosen to be non-useless, and all (except
$y$ and $z$) will avoid each $M_j$. Remove this copy of~\graph. Then
also greedily take out $r-3$ further disjoint copies of \graph,
avoiding the $M_j$ and all useless vertices, in such a way that
altogether they respect the proportions of the~$A_i$. Note that we
can find these copies greedily since the $(q+1)$-partite graph
induced by the $A_i$ is almost complete. We continue doing this
until no exceptional vertices are left in $A_{r-1}$. The fact that
$M_i$ has $s_i$ edges ensures that we will always have an edge left
in the appropriate matching for each exceptional vertex
in~$A_{r-1}$.

Now for all other exceptional vertices, proceed using the argument for
the case when $q \leq r-3$. In this way we will remove all the exceptional
vertices.

\medskip

\noindent So in both cases we will obtain sets $A_1', \ldots
,A_{q+1}'$ not containing any exceptional vertices. We now want to
remove any remaining useless vertices. Before dealing with the
exceptional vertices, each useless but non-exceptional vertex in
$A_i$ had at least $\tau^{1/3}|A_j|/2$ neighbours in $A_j$ for each
$j \neq i$. Also, we had at most $\tau^{2/3} n$ useless vertices,
and therefore also at most this many exceptional vertices. So we
have taken out at most $r(r-2)\tau^{2/3}n$ vertices. Thus each
remaining vertex $x \in A_i'$ still has at least
$\tau^{1/3}|A_j'|/3$ neighbours in $A_j'$ for each $j \neq i$, which
is much larger than the number of $j$-useless vertices.

Ideally, for a useless vertex $x\in A_i'$ we would like to pick
neighbours in each other class greedily so that together these
vertices form a copy of $K_r^-$ with, say, two vertices in $A_1'$,
$r-q-1$ vertices in $A_{q+1}'$ and one vertex in each other $A_j'$.
The problem is that the neighbours of~$x$ may avoid
a substantial proportion of $A_{q+1}'$, and so in particular may not
include any of the copies of $K_{r-q-1}$ which we know are contained
in $A_{q+1}$ (and therefore in~$A_{q+1}'$).

So instead, we proceed as follows. We first deal with all the
vertices which have too few neighbours in $A_{q+1}'$. Let $U$ be the
set of vertices in $A_1' \cup \ldots \cup A_q'$ which originally had
at most $(1-\tau^{1/4})|A_{q+1}|$ neighbours in $A_{q+1}$. In
particular, all these vertices are useless. Note that a vertex $x
\in U \cap A_i'$ (where $i \leq q$) still has at least
$\tau^{1/3}|A_i'|/3$ neighbours in $A_i'$. For each such vertex $x$ in
turn we proceed as follows. We first move~$x$ into $A_{q+1}'$.
Then we will greedily find a copy of $K_r^-$ which avoids~$x$ and meets
each $A_j'$ with $j \leq q$ in precisely one vertex. Note that similarly
as in~(\ref{eqmindegAq+1}) one can show that
\begin{equation}\label{eqmindegA'q+1}
\delta(G[A_{q+1}'])\geq
\left(1-\frac{1}{r-q-2}+\frac{c(r)}{2}\right)|A_{q+1}'|.
\end{equation}
So we may apply the Erd\H{o}s-Stone theorem to find the necessary
copy of $K_{r-q}^-$ in $A_{q+1}'$ avoiding~$x$ as well as all the
$(q+1)$-useless vertices. We can extend it to the desired copy of
$K_r^-$, also avoiding all the
useless vertices. Remove this copy of $K_r^-$. In effect, we have
removed two vertices from $A_i'$ (one vertex in the copy of~$K_r^-$ and $x$),
$r-q-1$ vertices from $A_{q+1}'$ and one vertex
from each other~$A_j'$. We can also find $r-3$ further disjoint
copies of $K_r^-$ in such a way that altogether these copies respect
the proportions of the $A_i'$. Remove these copies.
Repeating this for each vertex $x \in U$, in total we move or remove
at most $\tau^{1/2}n$ vertices. We denote~by $A_i''$ the sets thus
obtained from the~$A_i'$.

The effect of moving the vertices of $U$ and taking out these copies
of $K_r^-$ is that all vertices (except those in
$A_{q+1}''$) are joined to almost all of $A_{q+1}''$. The
vertices in~$U$ may now be $(q+1)$-useless, but are certainly
non-exceptional.

Now consider any useless vertex $x \in A_i''$ where $i \neq q+1$.
Let $A_j''$ be the vertex set in which $x$ has the lowest number
of neighbours, not including $j=i,q+1$. (Note that such a $j$ exists
since if $q=1$, a useless vertex $x \in A_1$ would have been in $U$,
so we would already have dealt with it.) Pick non-useless neighbours
$y$ and $z$ of $x$ in $A_j''$. (Such neighbours exist since $x$ was
not $j$-exceptional.) Recall that each of~$x$, $y$ and~$z$ is joined
to almost all of $A_{q+1}''$. Since $A_{q+1}''$ is almost as large as $A_{q+1}$
it follows that many of the copies of $K_{r-q-1}$ chosen
after~(\ref{eqmindegAq+1}) lie in the common neighbourhood
of $x$, $y$ and $z$, and so form a copy of $K_{r-q+2}^-$ together
with $x$, $y$ and $z$. Pick such a copy. Now note that the choice
of~$j$ implies that~$x$ is joined to at least $|A_{\ell}''|/3$
vertices in $A_{\ell}''$ for each $\ell\neq i,j,q+1$. So we can
greedily extend this copy of $K_{r-q+2}^-$ to a copy of $K_r^-$ in
$G$ by picking one non-useless vertex in every other~$A_{\ell}''$.
%
%
We then greedily find $r-3$ further disjoint copies of $K_r^-$
avoiding all the useless vertices so that together with the copy
just found, these copies of $K_r^-$ respect the proportions of the
$A_i''$. Remove all these copies of~$K_r^-$.

For a $(q+1)$-useless vertex~$x$, we perform a similar process, except
that $x$ is already in $A_{q+1}''$, so we find non-useless
neighbours $y$ and $z$ of $x$ in $A_j''$ and find a copy of
$K_{r-q-1}$ in $A_{q+1}''$ which contains $x$ and lies in the common
neighbourhood of $y$ and $z$. We can do this since~(\ref{eqmindegA'q+1}) implies that
$$
\delta(G[A_{q+1}''])\geq \left(1-\frac{1}{r-q-2}+\frac{c(r)}{3}\right)|A_{q+1}''|.
$$
(Note that in particular this bound applies to the degree of $x$ in $A_{q+1}''$.) So we
can successively pick common non-useless neighbours of $x$, $y$ and $z$ in $A_{q+1}''$
to construct the necessary $K_{r-q-1}$ containing $x$. Together with $y$
and $z$ this forms a copy of $K_{r-q+1}^-$ which we extend suitably
to a copy of $K_r^-$. As before we then find further disjoint copies of
$K_r^-$ such that together all these copies respect the proportions of the $A_i''$.
We can repeat this process until no useless vertices are left. The
fact that there are not too many useless vertices will ensure that
all our calculations remain valid.

Finally, if $k>0$, we remove $k$ further disjoint copies of $K_r^-$ to ensure
that the sets $A_1^*,\ldots ,A_{q+1}^*$ thus obtained from the
$A_i''$ are $(q,n^*)$-canonical where $n^*:=|A_1^*\cup\ldots\cup A_{q+1}^*|$.
This can be done because of our modification
of the $A_i$ at the beginning of the proof. Since the $A_i^*$
contain neither exceptional nor useless vertices and since we have
not removed too many vertices, it is easy to check that the $A_i^*$
satisfy all the conditions of the lemma. \qed

\section{Proof of Theorem \ref{Klminusthm}}\label{sec:end}
Recall that $B^*=B^*(K_r^-)$ denotes the bottle graph of~$K_r^-$.
Fix constants
$0<\tau_1 \ll \tau_2 \ll \ldots \ll \tau_{r-1} \ll 1/r$. Let $G$ be
the graph given in Theorem~\ref{the:main}. Let
$q\le r-2$ be maximal such that the conditions of Lemma~\ref{lem:tidy} are
satisfied with $\tau:=\tau_q$. As already observed
in Section~\ref{sec:overview}, by Theorem~\ref{the:nonex} we may assume
that $q\ge 1$. To prove Theorem~\ref{Klminusthm}, we apply first
Lemma~\ref{lem:tidy} with this choice of $q$ to obtain a subgraph
$G^*$ of $G$ and a $(q,|G^*|)$-canonical partition
$A_1^*,\ldots,A_{q+1}^*$ of~$V(G^*)$. Our definition of~$q$ will ensure that
if $q\neq r-3$ then the graph induced by~$A^*_{q+1}$ does not look
like one of the extremal graphs and so we can apply Theorem~\ref{the:nonex}
to it in order to find a perfect $B_1$-packing, where $B_1$ is the
spanning subgraph of $B^*_1$ defined below. (Recall that $B^*_1$
is the $(r-q-1)$-partite subgraph of $B^*$ obtained by deleting $q$
of the large vertex classes.) In the case when $q=r-3$ the graph $G^*[A^*_{q+1}]$
might violate condition~(ii) of Theorem~\ref{the:nonex}. So in this case
we will apply Theorem~\ref{the:nonex} to the ``almost-components'' of
$G^*[A^*_{q+1}]$ instead.

Recall that $A_1^*,\ldots,A_q^*$ all have the same size, which is
a multiple of $r-1$ (the size of a large class of the bottle
graph~$B^*$). The size of $A_{q+1}^*$ is a multiple of~$|B^*_1|$.
Our aim is to find a perfect $B_1$-packing in $G^*[A_{q+1}^*]$, where
$B_1$ is the graph consisting of $q$ vertex disjoint copies of
$K_{r-q-1}$ together with $r-q-2$ vertex disjoint copies of
$K_{r-q}^-$. We think of these copies as being arranged into an
$(r-q-1)$-partite graph with one vertex set of size $r-2$ and
$r-q-2$ vertex sets of size $r-1$. Thus $B_1\subseteq B^*_1$ and
the vertex classes of $B_1$ have the same sizes as those
of~$B^*_1$. This $B_1$-packing in $G^*[A_{q+1}^*]$ will then be
extended to a perfect \graph-packing in~$G^*$.

\begin{lemma}\label{lem:B1pack} We can take out from $G^*$ at most
$\tau^{1/3}n^*$ disjoint copies of $K_r^-$ to obtain subsets
$A^\diamond_1,\dots,A_{q+1}^\diamond$ of $A^*_1,\dots,A_{q+1}^*$ and
a subgraph $G^\diamond$ of $G^*$ such that the sets
$A^\diamond_1,\dots,A_{q+1}^\diamond$ are
$(q,|G^\diamond|)$-canonical and such that
$G^\diamond[A_{q+1}^\diamond]$ contains a perfect $B_1$-packing.
\end{lemma}
\proof Note that in the case when $q=r-2$ the graph $B_1$ just consists of $r-2$
isolated vertices, and the existence of a perfect $B_1$-packing is trivial
since $r-2$ divides $|A^*_{r-1}|$. In the case when $q\le r-3$ the
proof of Lemma~\ref{lem:B1pack} will invoke the non-extremal result,
Theorem~\ref{the:nonex}, with $\tau_{q+1}$ playing the role of $\tau_0$
there. It is for this reason that we will need the term~$-\theta n$
in the minimum degree condition in Theorem~\ref{the:nonex}.
Finally, note that hcf$(B_1)=1$ (even in the case when
$B_1$ is bipartite, i.e.~when $q=r-3$).%
     \COMMENT{$B_1$
has a class of size $r-2$ and a class of size $r-1$, so hcf$_\chi
(B_1) = 1$. For the bipartite case, note also that since $r \geq 4$,
there is at least one $K_{r-q-1}$ and at least one $K_{r-q}^-$ in
$B_1$, so hcf$_c(B_1)=1$.}
Let $s:=r-q-1\ge 2$. Thus $B_1$ is an
$s$-partite graph. Observe that $\chi_{cr}(B_1)=\chi_{cr}(B^*_1)=\frac{s(r-1)-1}{r-1}$.%
     \COMMENT{$|B_1|=(r-q-1)(r-1)-1$, $\sigma(B_1)=r-2$, $\chi(B_1)=r-q-1$.\\
Thus $\chi_{cr}(B_1)=\frac{(r-q-2)(s(r-1)-1)}{s(r-1)-1-(r-2)}=
\frac{(s-1)(s(r-1)-1)}{(s-1)(r-1)}=\frac{s(r-1)-1}{r-1}$}
Using~(i) and~(ii) of Lemma~\ref{lem:tidy}, similarly as in~(\ref{eqmindeg1})
and the first inequality in~(\ref{eqmindeg2})
one can show that
\begin{equation}\label{mindeg1}
\delta(G[A_{q+1}^*])\geq \left(1-\frac{1}{\chi_{cr}(B_1)}-
\tau_q^{1/4}\right)|A_{q+1}^*|=
\left(\frac{(s-1)(r-1)-1}{s(r-1)-1} -\tau_q^{1/4}\right)
|A_{q+1}^*|.
\end{equation}
So the minimum degree condition of
Theorem~\ref{the:nonex} is satisfied with $\theta:= \tau_q^{1/4}\ll
\tau_{q+1}$. Our choice of~$q$ implies that $G^*[A^*_{q+1}]$ satisfies
condition~(i) of Theorem~\ref{the:nonex} (with $\tau_0:=\tau_{q+1}$).
Thus in the case when $s>2$ we can apply Theorem~\ref{the:nonex} to find
a perfect $B_1$-packing in $G^*[A^*_{q+1}]$.

So we only need to consider the case when~$s=2$.
In this case $B_1$ is the bipartite graph consisting of $r-3$ disjoint edges and
one path of length~2, and we are done if condition~(ii)
of Theorem~\ref{the:nonex} holds. So suppose not and we do have some set
$C_1 \subseteq A_{q+1}^*$ with $d(C_1,A_{q+1}^*\setminus C_1) \leq \tau_{q+1}$.
Define $C_2:= A_{q+1}^* \setminus C_1$. Then there is a vertex $x \in C_1$
which has at most $\tau_{q+1}|C_2| \le \tau_{q+1}|A_{q+1}^*|$ neighbours in~$C_2$.
Together with~(\ref{mindeg1}) this shows that
$|C_1|>\delta(G^*[A^*_{q+1}])-\tau_{q+1}|A^*_{q+1}| \ge |A_{q+1}^*|/3$.
Similarly, $|C_2|>|A_{q+1}^*|/3$.

We now aim to show that by moving a few vertices,
we can achieve that each vertex in $C_1$ has few
neighbours in $C_2$ and vice versa. (This in turn will imply that
the graphs induced by both~$C_1$ and $C_2$ have large minimum degree.)
Call a vertex $x\in C_i$ in \emph{useless} if it has at most
$|C_i|/3$ neighbours in~$C_i$. By~(\ref{mindeg1}) every
such $x$ has at least $|C_j|/3$ neighbours in the other class~$C_j$.
Furthermore, the low density between~$C_1$ and~$C_2$
shows that there are at most $\tau_{q+1}^{3/4}|A_{q+1}^*|$ useless
vertices. We move each useless vertex into the other class and still
denote the classes thus obtained by~$C_1$ and~$C_2$.
Then $d(C_1,C_2)\leq \tau^{2/3}_{q+1}$.%
     \COMMENT{We have moved at most $\tau_{q+1}^{3/4}|A_{q+1}^*|$ vertices,
     so we have altered at most $\tau_{q+1}^{3/4}|A_{q+1}^*|^2$ edges. Thus
     \begin{align*}
     d(C_1,C_2) & = \frac{e(C_1,C_2)}{|C_1||C_2|}\\
     & \leq \tau_{q+1} +
     \frac{\tau_{q+1}^{3/4}|A_{q+1}^*|^2}{|C_1||C_2|}\\
     & \leq \tau_{q+1} + 9 \tau_{q+1}^{3/4}\\
     & \leq \tau_{q+1}^{2/3}
     \end{align*}
     }
Now call a vertex $x$ in either class \emph{bad} if it has at
least a $\tau_{q+1}^{1/6}$-fraction of the vertices in the other
class as neighbours. Clearly there are at most
$\tau_{q+1}^{1/2}|A_{q+1}^*|$ bad vertices. For each bad
vertex~$x\in C_i$ in turn we greedily choose a copy of~$B_1$ in~$C_i$
containing~$x$ such that these copies are disjoint for distinct
bad vertices. (Use that $\delta(G^*[C_i])\ge |C_i|/4$ for $i=1,2$
and the fact that~$B_1$ consists only of edges and a path of
length~2 to see that such copies can be found.) By removing these
copies of $B_1$, we end up with two sets $C_1'$ and $C_2'$ which
do not contain bad vertices. So each vertex in~$C'_1$ has at most
$2\tau_{q+1}^{1/6}|C_2'|$ neighbours in~$C_2'$ and vice versa.
Since $|C'_i|\ge |A^*_{q+1}|/4$ for $i=1,2$ (and thus also $|C_i'|
\leq 3|A_{q+1}^*|/4$ for $i=1,2$) this in turn implies that
\begin{equation}\label{eqminCs}
\delta(G^*[C_i'])\stackrel{(\ref{mindeg1})}{\ge}
\left(1-\frac{1}{\chi_{cr}(B_1)}-\tau_{q+1}^{1/7}\right)\frac{4|C_i'|}{3}
>\left(1-\frac{1}{\chi_{cr}(B_1)}\right)|C_i'|.
\end{equation}

We now aim to take out a few further copies of $K_r^-$ from $G^*$ to
ensure that both~$|C'_1|$ and~$|C'_2|$ are divisible by $|B_1|$.
As observed at the beginning of this section, $|A^*_{q+1}|$ is
divisible by~$|B_1|$. Thus $|C_1'|+|C_2'|$ is also divisible
by~$|B_1|$. Assume first that $|C_1'|=m|B_1|-1$ for some $m\in\N$.
We aim to remove
$2(r-2)$ disjoint copies of $K_r^-$ from $G^*$ in such a way that we
remove $2(r-1)$ vertices from every $A^*_i$ with $i\le r-3$,
$(r-1)+(r-2)-1$ vertices from $C'_1$ and $(r-1)+(r-2)+1$ vertices
from~$C'_2$. Then the sizes of the remaining subsets of~$C'_1$ and~$C'_2$
will be divisible by~$|B_1|$. Moreover, since the $A^*_i$ were $(q,|G^*|)$-canonical,
and since altogether we remove $2((r-1)+(r-2))$ vertices from
$A^*_{q+1}$, the remaining subsets will still induce a canonical
partition of the remaining subgraph of~$G^*$.

The way we remove the above copies of $K_r^-$ is as follows:
Greedily find $r-2$ disjoint copies of $K_r^-$ with two vertices in $C_1'$,
two vertices in $A_i^*$ and one
vertex in each $A_j^*$ with $1 \leq j \leq r-3$ and $j\neq i$. For
each of these copies of $K_r^-$ the index~$i$ will be different except
that $i=1$ will be chosen twice. Also find $r-4$ disjoint copies of $K_r^-$ with two
vertices in $C_2'$, two vertices in $A_i^*$ and one vertex in each
$A_j^*$ with $1 \leq j \leq r-3$ and $j\neq i$.  The choices of $i$
will be between~2 and $r-3$, and no~$i$ will be chosen twice.
Finally, find two copies of $K_r^-$ with three vertices in $C_2'$ and
one in each $A_i^*$ for $1 \leq i \leq r-3$.

In the general case (i.e.~when $|C'_i|\equiv t\mod |B_1|$), we
simply repeat this procedure $t$ times to even out the residues
modulo $|B_1|$ between $|C_1'|$ and $|C_2'|$. We denote the
remaining subsets by $A^\diamond_i$ and $C^\diamond_i$ and the
remaining subgraph by~$G^\diamond$. We only need to perform the
above procedure at most $|B_1|-1$ times, so we are taking out a
bounded number of copies of $K_r^-$, which will not affect any of
the vertex degrees significantly. Thus each
$G^\diamond[C^\diamond_i]$ satisfies the minimum degree condition
in Theorem~\ref{the:nonex}. Indeed, the first inequality in~(\ref{eqminCs}) shows that
\begin{equation}\label{checkconds}
\delta(G^\diamond[C_i^\diamond]) \geq
\left(1-\frac{1}{\chi_{cr}(B_1)}-\tau_{q+1}^{1/8}\right)\frac{4|C_i^\diamond|}{3} \geq
\left(\frac{2}{5}-\tau_{q+1}^{1/8}\right)\frac{4|C_i^\diamond|}{3}\ge
\frac{51}{100}|C_i^\diamond|.
\end{equation}
This bound on the minimum degree also shows that each $C_i^\diamond$
cannot contain an almost independent set of size $|C_i^\diamond|/2$,
so condition~(i) of Theorem~\ref{the:nonex} is satisfied with room
to spare. To see that condition~(ii) also holds, observe that if
$C_i^\diamond$ is partitioned into $S_1$ and $S_2$, where
$0<|S_1|\leq |C_i^\diamond|/2 \leq |S_2|$, then the neighbours of
any vertex in $S_1$ cover a significant proportion (at least $1/50$)
of $S_2$, and so $d(S_1,S_2) \geq 1/50$. So condition~(ii) is
satisfied too. Thus we can apply Theorem~\ref{the:nonex} to each of
the subgraphs of $G^\diamond$ induced by $C^\diamond_1$ and
$C^\diamond_2$ to find perfect $B_1$-packings in
$G^\diamond[C^\diamond_1]$ and $G^\diamond[C^\diamond_2]$. Adding
back into $A_{q+1}^\diamond$ the vertices in the copies of $B_1$
which were removed when dealing with the bad vertices (and letting
$G^\diamond$ denote the subgraph of $G$ induced by the modified
$A_i^\diamond$), we still have a perfect $B_1$-packing in
$G^\diamond[A_{q+1}^\diamond]$, and $G-G^\diamond$ consists of those
copies of $K_r^-$ which we removed. Thus $G^\diamond$ and the
$A_i^\diamond$ are as required in the lemma.
\endproof

Our aim now is to extend the perfect $B_1$-packing in $G^\diamond[A^\diamond_{q+1}]$
to a perfect $K^-_r$-packing in~$G^\diamond$.%
     \COMMENT{Can't write $B^*$-packing in~$G^\diamond$ since we only
have a $B_1$-packing and not a $B^*_1$-packing.} To do this, we
define a ($q+1$)-partite auxiliary graph~$J$, whose vertices are the
vertices in $A_i^\diamond$ for all $1 \leq i \leq q$ together with
all the copies of~$B_1$ in the perfect $B_1$-packing of
$G^\diamond[A_{q+1}^\diamond]$. There will be an edge between
vertices from the $A_i^\diamond$'s whenever there was one in $G$,
and a vertex $x\in A_i^\diamond$ for $1 \leq i \leq q$ will be
joined to a copy of~$B_1$ whenever $x$ was joined to all the
vertices of this copy in $G$.

Let $H_{q,r-1}$ denote the complete
$(q+1)$-partite graph with $q$ classes of size $r-1$ and one
class of size~1.
We wish to find a perfect $H_{q,r-1}$-packing in~$J$. It is easy to see that
this then yields a perfect \graph-packing in~$G^\diamond$
and thus, together with all the copies of $K_r^-$ chosen earlier, a
perfect \graph-packing in~$G$.

The existence of such a perfect $H_{q,r-1}$-packing follows
immediately from Proposition~\ref{almostcomplete}. To see that we
can apply this proposition, note that Lemma~\ref{lem:tidy}(iii)
implies that in $G^*$ each vertex is adjacent to almost all vertices
in the other vertex classes and this remains true in~$G^\diamond$
since we only deleted a small proportion of the vertices after
applying Lemma~\ref{lem:tidy}. It follows immediately that every
vertex in~$J$ is adjacent to almost all vertices in the other vertex
classes of~$J$. Note also that the vertex classes of~$J$ have the
correct sizes since the sets $A^\diamond_1,\dots,A^\diamond_{q+1}$
are $(q,|G^\diamond|)$-canonical. This completes the proof of
Theorem~\ref{Klminusthm}.

\discard
\bigskip

The ideas in this proof can easily be adapted to a larger class of
graphs which we will now define. Call a colouring $c$ of $H$
\emph{appropriate} if $c$ is optimal and has a colour class of size
$\sigma(H)$. Given an appropriate colouring~$c$ of a graph~$H$ of
chromatic number~$\ell\ge 3$, let $x_1 \leq x_2 \leq \ldots \leq
x_\ell$ be the sizes of the colour classes. Define
$\mathcal{D}^{deg}(c)=\lbrace x_{i+1}-x_i \mid i=2,\ldots,\ell-1
\rbrace$. Let $\mathcal{D}^{deg}(H)$ be the union of all the sets
$\mathcal{D}(c)$ over all appropriate colourings~$c$ of~$H$. We
define $hcf_\chi^{deg}(H)$ to be the highest common factor of the
elements of~$\mathcal{D}^{deg}(H)$ (or $hcf_\chi^{deg}(H):=\infty$
if $\mathcal{D}^{deg}(H) = \lbrace 0 \rbrace$). Note that if
$hcf_\chi^{deg}(H)=1$ then
$hcf_\chi(H)=1$.%
     \COMMENT{What about bipartite graphs?}

Also define, instead of $B^*(H)$, a larger graph $D^*(H)$. Where the
vertex classes of $B^*$ were formed by taking the union over cyclic
permutations of the large classes of $H$, the vertex classes of
$D^*$ will formed by taking the union over all appropriate
colourings. Note that this automatically covers the cyclings of
large classes, since such a cycling would result in another
appropriate colouring. Thus, if the number of appropriate colourings
is $m$, then $\ell-1$ divides $m$ and $D^*(H)$ will have $\ell-1$
classes of size $m(|H|-\sigma)/(\ell-1)$ and one class of size
$m\sigma$.

For an appropriate colouring $c$ of $H$, define $H_q^c$ to be the
subgraph of $H$ induced by the $\ell-q$ smallest classes of $H$ in
this colouring. For $q \leq \ell-2$, define $D_q^*(H)$ to be the
subgraph of $D^*(H)$ induced by the small class and $\ell-q-1$ of
the large classes. Define $D_q(H)$ to be the analogous (to $B_1$)
non-complete graph as formed by taking the union of the $H_q^c$ over
all appropriate colourings $c$ of $H$.

Now suppose that $H$ is an $\ell$-partite graph, where $\ell \geq
3$, which satisfies the following conditions:

\medskip

\begin{itemize}
\item[{\rm (i)}] $hcf_\chi^{deg}(H)=1$.
\item[{\rm (ii)}] $H$ has an appropriate colouring in which some vertex,
$y_0$, lies in a class of size $\sigma$ and has exactly one
neighbour, $z_0$, in some other class.
\item[{\rm (iii)}] $H$ has an optimal colouring with a vertex class of size
$\sigma + 1$.
\item[{\rm (iv)}] hcf$(D_{\ell-2}(H))=1$.
\item[{\rm (v)}] hcf  $\lbrace |H_{\ell-2}^c| \mid$
$c$ is an appropriate colouring$\rbrace = 1$.
\end{itemize}

In this case the ideas in the proof of Theorem~\ref{the:main} also
show that $G$ has a perfect $H$-packing provided that the minimum
degree of $G$ satisfies
$$ \delta(G) \geq \left(1-\frac{1}{\chi_{cr}(H)}\right)|G|$$
and that $|G|$ is divisible by~$|H|$ and sufficiently large.%
\COMMENT{Conditions (ii)-(v) can be shown to be necessary. Given
(v), (iv) amounts to saying hcf$_\chi(D_{\ell-2}(H))\leq 2$} Details
will be given in~\cite{MPhil}.
\enddiscard

{\footnotesize
\bigskip
\noindent
Oliver Cooley, Daniela K\"uhn \& Deryk Osthus\\
School of Mathematics\\
Birmingham University\\
Edgbaston\\
Birmingham B15 2TT\\
UK\\
{\it E-mail addresses}: {\tt
\{cooleyo,kuehn,osthus\}@maths.bham.ac.uk} }

\end{document}